\documentclass[12pt]{amsart}
\usepackage{graphicx, overpic}
\usepackage[below]{placeins}
\usepackage[colorlinks=true, linkcolor=blue, citecolor=blue]{hyperref}
\usepackage[]{algorithm2e}

\usepackage[T1]{fontenc}

\usepackage{amsmath,amsthm,amscd,amssymb,eucal}

\usepackage{enumerate, amsfonts, latexsym, color, url}
\usepackage{epstopdf}

\usepackage{pinlabel}
\usepackage{tikz}

\makeatletter
\newsavebox{\@brx}
\newcommand{\llangle}[1][]{\savebox{\@brx}{\(\m@th{#1\langle}\)}%
  \mathopen{\copy\@brx\kern-0.5\wd\@brx\usebox{\@brx}}}
\newcommand{\rrangle}[1][]{\savebox{\@brx}{\(\m@th{#1\rangle}\)}%
  \mathclose{\copy\@brx\kern-0.5\wd\@brx\usebox{\@brx}}}
\makeatother

\begin{document}

\newtheorem{theorem}{Theorem}[section]
\newtheorem{lemma}[theorem]{Lemma}
\newtheorem{proposition}[theorem]{Proposition}
\newtheorem{corollary}[theorem]{Corollary}
\newtheorem{conjecture}[theorem]{Conjecture}
\newtheorem{question}[theorem]{Question}
\newtheorem{problem}[theorem]{Problem}
\newtheorem*{claim}{Claim}
\newtheorem*{criterion}{Criterion}

\theoremstyle{definition}
\newtheorem{definition}[theorem]{Definition}
\newtheorem{construction}[theorem]{Construction}
\newtheorem{notation}[theorem]{Notation}
\newtheorem{object}[theorem]{Object}
\newtheorem{operation}[theorem]{Operation}

\theoremstyle{remark}
\newtheorem{remark}[theorem]{Remark}
\newtheorem{example}[theorem]{Example}

\numberwithin{equation}{subsection}

\newcommand\id{\textnormal{id}}

\newcommand\Z{\mathbb Z}
\newcommand\R{\mathbb R}
\newcommand\C{\mathbb C}
\newcommand\CC{\mathbf C}
\newcommand\BB{\mathbf B}
\newcommand\TT{\mathbf T}
\newcommand\PP{\mathcal P}
\newcommand\W{\mathcal W}
\newcommand\RR{\mathcal R}
\newcommand\Aut{\textnormal{Aut}}

\title{Wiggle Island}

\author{Danny Calegari}
\address{University of Chicago \\ Chicago, Ill 60637 USA}
\email{dannyc@math.uchicago.edu}
\date{\today}

\begin{abstract}
A {\em wiggle} is an embedded curve in the plane that is the attractor of an
iterated function system associated to a complex parameter $z$. We show the
space of wiggles is disconnected --- i.e.\/ there is a {\em wiggle island}.
\end{abstract}

\maketitle
\setcounter{tocdepth}{1}

\section{Squiggles and wiggles}

A {\em squiggle} is a continuous map $w_z:I \to \C$ where $I:=[0,1]$,
depending on a complex parameter $z$ with $|z|<1$ and $|1-z|<1$ in the following way. 

\begin{figure}[hpbt]
\centering
\includegraphics[scale=0.35]{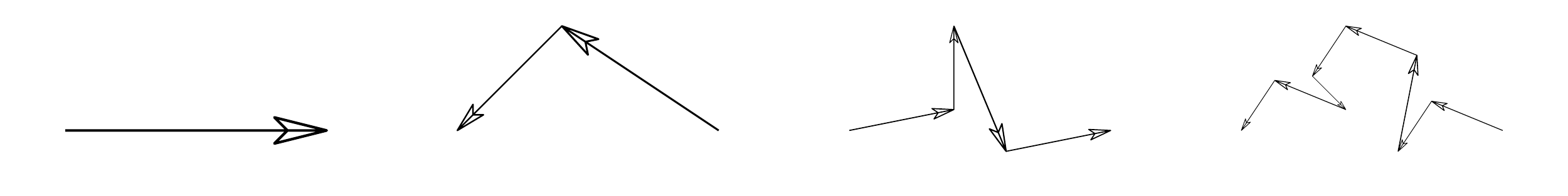}
\caption{$w^j_z$ for $z=0.4+0.4i$ and $j=0,1,2,3$}\label{wiggle_iterates}
\end{figure}

If we define the two maps
$$f_z: x \to -zx+z, \quad g_z: x \to (z-1)x+1$$
and inductively define a sequence of maps $w^j_z:I \to \C$ by
\begin{enumerate}
\item{$w^0_z:[0,1] \to \C$ is the identity map; and}
\item{for $j\ge 1$ we first define $\widetilde{w}^j_z:I \to \C$ by
\[ \widetilde{w}^j_z(t) := \begin{cases}
g_z\bigl( w^{j-1}_z(2t)\bigr) & \text{for } 0\le t\le 1/2, \\
f_z\bigl( w^{j-1}_z(2t-1)\bigr) & \text{for } 1/2\le t\le 1,
\end{cases} \]
and then set $w^j_z(t):=\widetilde{w}^j_z(1-t)$;}
\end{enumerate}
then $w_z$ is the limit of $w^j_z$ as $j \to \infty$. This limit exists
because when $|z|<1$ and $|1-z|<1$ both $f_z$ and $g_z$ are uniformly
contracting. See Figure~\ref{wiggle_iterates}
illustrating $w^j_z$ (up to reverse of orientation) for $z=0.4+0.4i$ and $j=0,1,2,3$.

In words, $w^j_z$ is the map obtained by concatenating the two maps
$x \to g_z w^{j-1}_z x$ and $x \to f_z w^{j-1}_z x$ (in that order), and then
reversing the orientation (i.e.\/ precomposing with $x \to 1-x$).
Since $f_z(0)=z=g_z(1)$ the two cases in the definition of $\widetilde{w}^j_z$
agree at $t=1/2$, so $\widetilde{w}^j_z$ and therefore $w^j_z$ is continuous;
and since $f_z(1)=0$ and $g_z(0)=1$, induction on $j$ gives
$$w^j_z(0)=f_z\bigl(w^{j-1}_z(1)\bigr)=0, \quad w^j_z(1)=g_z\bigl(w^{j-1}_z(0)\bigr)=1$$
for every $j$. This last point is the purpose of the reversal: without it the
two endpoints would be interchanged at each stage, and the sequence of maps
would not converge.

\begin{figure}[hpbt]
\centering
\includegraphics[scale=0.3]{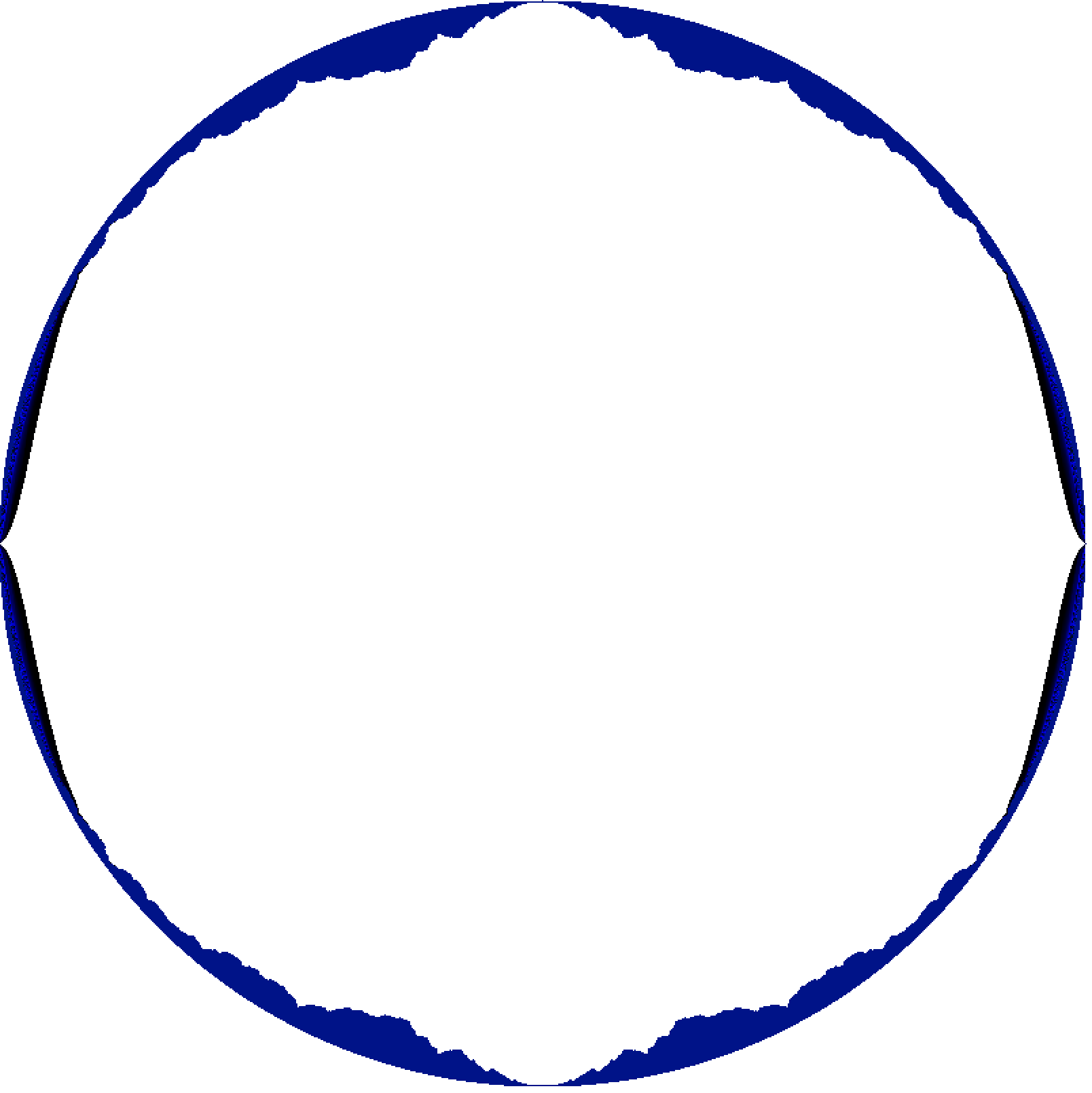}
\caption{$\W$ (in white) as a subset of $\lbrace z: |z-1/2|<1/2\rbrace$.}\label{W}
\end{figure}

A {\em wiggle} is a squiggle which is an embedding. The set of wiggles $w_z$
is parameterized by an open subset 
$\W:= \lbrace z \text{ such that } w_z \text{ is a wiggle}\rbrace$ of $\C$.
A wiggle with parameter $z$ has Hausdorff dimension $d$
where $|z|^d + |1-z|^d = 1$; thus $\W$ is a subset of the open disk of radius
$1/2$ centered at $1/2$ (because a squiggle might intersect itself
its Hausdorff dimension only satisfies the inequality 
$|z|^d + |1-z|^d \ge 1$).
Denote the complement of $\W$ in this disk by $\RR$.
Figure~\ref{W} depicts $\W$ (in white) as a subset of this disk (note how $\W$ very
nearly fills the entire disk!)

The set $\W$ contains one {\em big component}, the connected component of $1/2$;
this big component contains the real interval $(0,1)$ and the imaginary 
interval $(0.5-0.5i,0.5+0.5i)$. However, the big component is not all of $\W$:

\begin{theorem}[Wiggle Island]\label{theorem:wiggle_island}
$\W$ is not connected.
\end{theorem}

\begin{figure}[hpbt]
\centering
\includegraphics[scale=0.25]{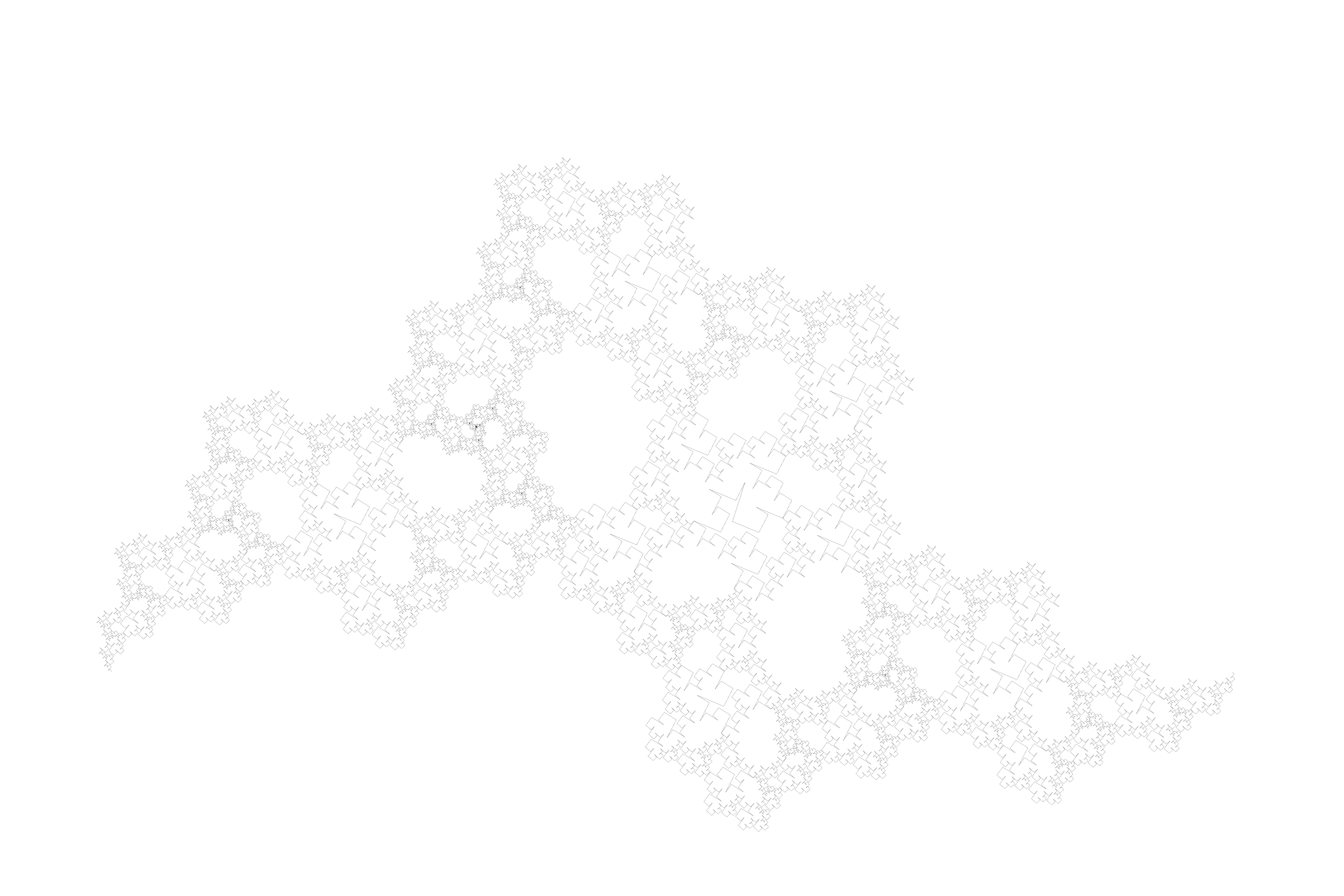}
\caption{An approximation to $w_z$ for $z=0.3409+0.43486i$; this is an embedded arc.}\label{wiggle_curve}
\end{figure}

Figure~\ref{wiggle_curve} depicts a wiggle $w_z$ for $z$ in an `island'
component of $\W$ centered at approximately $z=0.3409+0.43486i$. The island is
invisible at the resolution of Figure~\ref{W}; a zoomed in image of the island is
Figure~\ref{W_island_rigorous}.

The wiggle $w_z$ for $z$ as above {\em may not} be deformed through 
wiggles to $w_{1/2}$ (i.e.\/ the unit interval).

\begin{remark}[The Carpenter's Rule Problem]
The {\em Carpenter's Rule Problem}, first posed by Stephen Schanuel and George Bergman 
in the early 1970's, asks whether every embedded planar polygonal arc may be
{\em straightened} (i.e.\/ moved through embeddings to a straight arc by
changing the angles but not the lengths of the segments).
Connelly, Demaine and Rote showed in 2000 \cite{Connelly_Demaine_Rote} that
the answer to the Carpenter's Rule Problem is {\em yes}.

The {\em Wiggle Problem} asks analogously whether every wiggle may be straightened
through a family of wiggles; Theorem~\ref{theorem:wiggle_island} says that the
answer to the Wiggle Problem is {\em no}.
\end{remark}

\subsection{Two generator IFS}

Since $f_z$ and $g_z$ are contractions of $\C$, the pair $\lbrace f_z,g_z \rbrace$
is an {\em iterated function system} (IFS), and as such has a unique {\em attractor}:
the unique nonempty compact set $K_z \subset \C$ with $K_z = f_zK_z \cup g_zK_z$.
Directly from the construction of $w_z$ we have
$w_z(I) = f_zw_z(I) \cup g_zw_z(I)$, and $w_z(I)$ is nonempty and compact; so
$$K_z = w_z(I)$$
for every $z$ with $|z|<1$ and $|1-z|<1$. In other words a squiggle is nothing
other than the attractor of a two generator linear IFS, and $\W$ is the locus of
parameters for which this attractor is an embedded arc.

It is worth saying which two generator systems arise in this way. Let $\lambda,\mu$
be complex numbers of modulus less than $1$ and put $f(x) := \lambda(x-1)$ and
$g(x) := \mu x + 1$, so that $f(1)=0$ and $g(0)=1$ and the attractor $K$ contains
$0$ and $1$. Then the two halves $fK$ and $gK$ of $K$ share an endpoint --- i.e.\/
$f(0)=g(1)$ --- if and only if $\lambda + \mu = -1$; and taking $\lambda = -z$
recovers $f=f_z$ and $g=g_z$. Thus our family is precisely the $1$ complex
dimensional slice $\lambda+\mu=-1$ of the $2$ complex dimensional family of pairs
of linear contractions of $\C$: the slice on which the attractor is a chain of two
links, and therefore has a chance to be an arc.

Attractors of two map affine systems obeying this chain condition are classical.
They were studied by de Rham \cite{de_Rham_1956,de_Rham_1957}, who obtained them as
limits of a corner cutting procedure and was principally interested in the
parameters for which the limit is smooth; curves in this family are often called
{\em de Rham curves}. The general question of when the attractor of such a system
is a Jordan arc is taken up by Aseev, Tetenov and Kravchenko
\cite{Aseev_Tetenov_Kravchenko}, who call such a system a {\em zipper} and record
the orientations of its pieces in a {\em signature}. In this language a squiggle is
the attractor of a two map zipper, the signature recording the reversals of
orientation that appear in the definition of $w^j_z$.

Two generator linear IFS have been studied intensively, and the parameter spaces
that arise are typically intricate. The best known example is the
{\em Barnsley--Harrington Mandelbrot set}
$$\mathcal{M} := \lbrace \lambda \text{ with } 0<|\lambda|<1 \text{ such that the
attractor of } \lbrace \lambda x - 1, \lambda x + 1 \rbrace \text{ is connected}
\rbrace,$$
introduced in \cite{Barnsley_Harrington}; here the two generators have the {\em same}
ratio $\lambda$. Bousch \cite{Bousch} proved that $\mathcal{M}$ is connected and
locally connected. Bandt \cite{Bandt} discovered that $\mathcal{M}$ has interior
points, and conjectured that the interior of $\mathcal{M}$ is dense in $\mathcal{M}$
away from the real axis; this
is the conjecture proved in \cite{Calegari_Koch_Walker}, by the method of {\em traps}
that we adapt below. Solomyak \cite{Solomyak} studied the asymptotic self-similarity
of $\mathcal{M}$ near its frontier. A parameter $\lambda$ lies in $\mathcal{M}$ if
and only if
$\sum_n a_n \lambda^n = 0$ for some sequence $a_n \in \lbrace 0,\pm 1\rbrace$ not
identically zero, so that $\mathcal{M}$ is closely related to the sets of zeros of
polynomials
with restricted coefficients studied by Odlyzko and Poonen \cite{Odlyzko_Poonen}.

Unpublished work of Astala, Rohde and Schramm studies exactly the family
$\lbrace f_z,g_z\rbrace$, and exactly the question of for which parameters the
attractor is embedded. In their normalization the arc joins $-1$ to $1$ and the
parameter is the position $p$ of the image of the midpoint of $I$; since
$z = w_z(1/2)$ the two normalizations are related by $p=2z-1$, which takes the disk
$|z-1/2|<1/2$ to the unit disk. They study the set
$$J:=\lbrace p \text{ such that the attractor is a simple arc}\rbrace,$$
which corresponds to $\W$ under $p=2z-1$. Their motivation is a conjecture of
Astala from 1994 \cite{Astala}, that
$$\sup \dim f(\R) = 1+k^2$$
where the supremum is taken over all $k$ quasiconformal $f$; equivalently, that
$$\sup_F \dim F([-1,1],k) = 1+k^2$$
where this supremum is taken over all {\em holomorphic motions} $F$ of $[-1,1]$.
Any conformal map $\varphi$ from the unit
disk into $J$ gives rise to a holomorphic motion $(x,\zeta) \to F(x,\varphi(\zeta))$
of $[-1,1]$, and therefore to a lower bound for Astala's conjecture; with the
assistance of a computer they show that $J$ contains the disk of radius $0.91$ about
$0$, and deduce the lower bound $1+0.59k^2$. See \cite{Astala_Rohde_Schramm}.

Theorem~\ref{theorem:wiggle_island} bears directly on this strategy. Since $\W$,
and therefore $J$, is disconnected, $J$ is not conformally a disk, and any
$\varphi$ as above must omit part of $J$. In fact Wiggle Island sits at
$p = 2z-1 \sim -0.3182+0.8697i$, of modulus $\sim 0.9261$; so it lies in the thin
annulus $0.91<|p|<1$ that is not reached by the disk above, and no round disk
about $0$ contained in $J$ can reach it.

\section{Certifying in and out}

In this section we give two {\em stable numerical} criteria to certify
that $z\in \W$ resp. $z\in \RR$. 

\subsection{Certifying $z\in \W$}\label{subsection:certifying_W}

Fix $z$ and write $\gamma:=w_z(I)$ and abbreviate $f_z$ and $g_z$ by $f$ and $g$
so that $\gamma = f\gamma \cup g\gamma$.
Let $S_n$ denote the set of words of length $n$ in the alphabet $\lbrace f,g\rbrace$
and let $S=\cup_n S_n$.
By abuse of notation we think of $u \in S_n$ as a map, obtained by 
composing $f$ or $g$ according to the letters of $u$. 
Let $fS_{n-1}$ resp. $gS_{n-1}$ represent words of length $n$ beginning with
$f$ and $g$ respectively; more generally juxtaposition denotes concatenation,
so that $fS_k$ is the set of words of length $k+1$ beginning with $f$.

Applying the relation $\gamma = f\gamma \cup g\gamma$ repeatedly we obtain, by
induction,
$$f\gamma = \bigcup_{u\in fS_{n-1}} u\gamma, \qquad
g\gamma = \bigcup_{v\in gS_{n-1}} v\gamma$$
for all $n>1$.

\begin{lemma}\label{lemma:z_in_W}
$z\in \W$ if and only if $f\gamma \cap g\gamma = \lbrace z\rbrace$.
\end{lemma}
\begin{proof}
The map $w_z$ takes $[0,1/2]$ onto $f\gamma$ and $[1/2,1]$ onto $g\gamma$, and
$w_z(1/2)=z$. Thus if $w_z$ is injective then
$$f\gamma \cap g\gamma = w_z([0,1/2]) \cap w_z([1/2,1])
= \lbrace w_z(1/2) \rbrace = \lbrace z \rbrace.$$

Conversely suppose $f\gamma \cap g\gamma = \lbrace z\rbrace$, and let
$A_0:=w_z^{-1}(0)$ and $A_1:=w_z^{-1}(1)$. We claim $A_0 = \lbrace 0 \rbrace$ and
$A_1 = \lbrace 1\rbrace$. Since $|z|<1$ and $|1-z|<1$ we have
$z\ne 0$ and $z\ne 1$; and $0=f(1)$ lies in $f\gamma$ while $1=g(0)$ lies in
$g\gamma$, so neither $0$ nor $1$ can lie in both pieces, and therefore
$A_0 \subset [0,1/2]$ and $A_1\subset [1/2,1]$. On $[0,1/2]$ we have
$w_z(s)=f(w_z(1-2s))$, and $f(x)=0$ only for $x=1$; so
$A_0 = \lbrace (1-b)/2 : b \in A_1\rbrace$. Likewise on $[1/2,1]$ we have
$w_z(t)=g(w_z(2-2t))$, and $g(x)=1$ only for $x=0$; so
$A_1 = \lbrace 1-a/2 : a\in A_0\rbrace$. Combining these,
$A_0 = \lbrace a/4 : a\in A_0\rbrace$, so that $a\in A_0$ implies $4^na \in A_0$ for
every $n$; since $A_0\subset [0,1/2]$ this forces $a=0$, and the claim follows.

Now suppose $w_z(s)=w_z(t)$ for some $s\ne t$. Recording for each parameter the
sequence of halves in which it lies at each stage of the subdivision, let $p$ be the
longest common prefix of the sequences associated to $s$ and $t$; then after
relabelling $w_z(s) \in pf\gamma$ and $w_z(t)\in pg\gamma$. Since $p$ is an injective
affine map, $p^{-1}w_z(s) \in f\gamma\cap g\gamma$ and therefore $w_z(s)=p(z)$. The
parameter separating the two subintervals also has image $p(z)$; as $s\ne t$, at
least one of $s,t$ differs from it, and applying $(pf)^{-1}$ or $(pg)^{-1}$ and using
$f(0)=z=g(1)$ this exhibits an element of $A_0$ other than $0$, or an element of
$A_1$ other than $1$, contrary to the claim. Hence $w_z$ is injective and
$z\in \W$.
\end{proof}

\begin{lemma}\label{lemma:ball}
Define 
$$R = \max\left( \frac {|z-1|}{2|1-|z||}, \frac {|z|} {2|1-|1-z||} \right)$$
and let $B$ be the ball of radius $R$ about $1/2$.
Then $\gamma$ is contained in $B$.
\end{lemma}
\begin{proof}
Both $f$ and $g$ take $B$ inside itself.
\end{proof}

Thus to show $z\in \W$ it would suffice to show, for any $n$, that 
$uB \cap vB = \emptyset$ for all $u\in fS_n$ and $v\in gS_n$ except for 
one specific pair for which $uB \cap vB = \lbrace z \rbrace$.
Unfortunately this is impossible, at least when $z$ is not real; necessarily
$z$ itself is in the interior of some $uB \cap vB$.

For $ffg(0)=z$ and $ggf(1)=z$; and if $z$ is not real then $|z|+|1-z|>1$, so that
$|z-1|>1-|z|$ and $|z|>1-|1-z|$. Both terms in the definition of $R$ therefore
exceed $1/2$; so $R>1/2$, and $0$ and $1$ lie in the interior of $B$. Since $ffg$
and $ggf$ are injective affine maps it follows that
$z=ffg(0) \in \textnormal{int}(ffgB)$ and $z=ggf(1) \in \textnormal{int}(ggfB)$, and
therefore $z \in \textnormal{int}(ffgB\cap ggfB)$, where $ffg\in fS_2$ and
$ggf\in gS_2$. (For real $z$ we have $z\in(0,1)$ and $R=1/2$ exactly, so that $0$
and $1$ lie on $\partial B$ and this argument does not apply; but there is nothing
to prove for such $z$, since then $\gamma=[0,1]$ and $(0,1)\subset \W$.)

The key observation is that $ffgf^{-1} = ggfg^{-1}$; in particular,
$ffg\gamma \cup ggf \gamma$ is a dilated copy of $\gamma$ itself.
Indeed the linear parts of $ffg$ and $f^{-1}$ are $z^2(z-1)$ and $-1/z$, while those
of $ggf$ and $g^{-1}$ are $-z(z-1)^2$ and $1/(z-1)$; so $ffgf^{-1}$ and $ggfg^{-1}$
both have linear part $z(1-z)$. Moreover $f^{-1}(z)=0$ and $g^{-1}(z)=1$, so by the
identities $ffg(0)=z$ and $ggf(1)=z$ above both maps fix $z$. Since an affine map of
$\C$ is determined by its linear part together with the image of one point, the two
maps agree; write $h$ for their common value, so that
$$h(x)-z = z(1-z)(x-z).$$
Thus $h$ is the similarity fixing $z$ with ratio $z(1-z)$, and it is a contraction
because $|z|<1$ and $|1-z|<1$. Now $hf = ffg$ and $hg=ggf$, so
$\gamma = f\gamma \cup g\gamma$ gives
$$h\gamma = ffg\gamma \cup ggf\gamma,$$
which is the dilated copy referred to above; and since $h$ is injective,
$$ffg\gamma \cap ggf\gamma = h(f\gamma\cap g\gamma).$$
Thus if $ffg\gamma \cap ggf \gamma$ contains a point other than $z$, then
$f\gamma \cap g\gamma$ contains a point other than $z$ which is not in
$ffg\gamma \cap ggf \gamma$.
For suppose $p \in ffg\gamma \cap ggf\gamma$ with $p\ne z$, and put
$p_n := h^{-n}(p)$. By the identity above $p_1 \in f\gamma\cap g\gamma$; and whenever
$p_n$ lies in $ffg\gamma\cap ggf\gamma$ the same identity puts $p_{n+1}$ in
$f\gamma \cap g\gamma$, so we may continue. No $p_n$ is $z$, since $h$ is injective
and fixes $z$; but
$$|p_n - z| = |z(1-z)|^{-n}|p-z|$$
tends to infinity, whereas $f\gamma\cap g\gamma$ is bounded. So the process must
stop, and where it stops is a point of $f\gamma \cap g\gamma$ other than $z$ which is
not in $ffg\gamma \cap ggf\gamma$.
We therefore obtain the following algorithm which,
if it terminates, certifies that $z$ is in the interior of $\W$:

\begin{algorithm}
Initialize $L$ to the set of pairs $(u,v) \in fS_2 \times gS_2 - (ffg,ggf)$

\While{$L\ne \emptyset$}{
\ForAll{$(u,v) \in L$}{

\eIf{$uB \cap vB = \emptyset$}{
remove $(u,v)$ from $L$
}{
replace $(u,v)$ with $\lbrace (uf,vf),(uf,vg),(ug,vf),(ug,vg) \rbrace$
}
}
}
\end{algorithm}

\subsection{Roots and $\RR$}

There is a second description of $\RR$, as the closure of a set of roots of
polynomials, which brings out the analogy with the Barnsley--Harrington set
$\mathcal{M}$ of \S~1.1.

For a word $u\in S$ let $p_u(z)$ be the value at the point $z$ of the affine map
$u$; since the coefficients of $f$ and $g$ are polynomials in $z$, so is $p_u$, of
degree $|u|+1$. These polynomials are generated from the seed $p_\emptyset(z)=z$ by
$$p_{fu}(z) = z - zp_u(z), \qquad p_{gu}(z) = 1+(z-1)p_u(z);$$
thus $p_f = z-z^2$, $p_g = z^2-z+1$, $p_{ff} = z^3-z^2+z$,
$p_{gf} = -z^3+2z^2-z+1$, and so on.

\begin{figure}[hpbt]
\centering
\begin{minipage}[c]{0.48\textwidth}
\centering
\includegraphics[width=\textwidth]{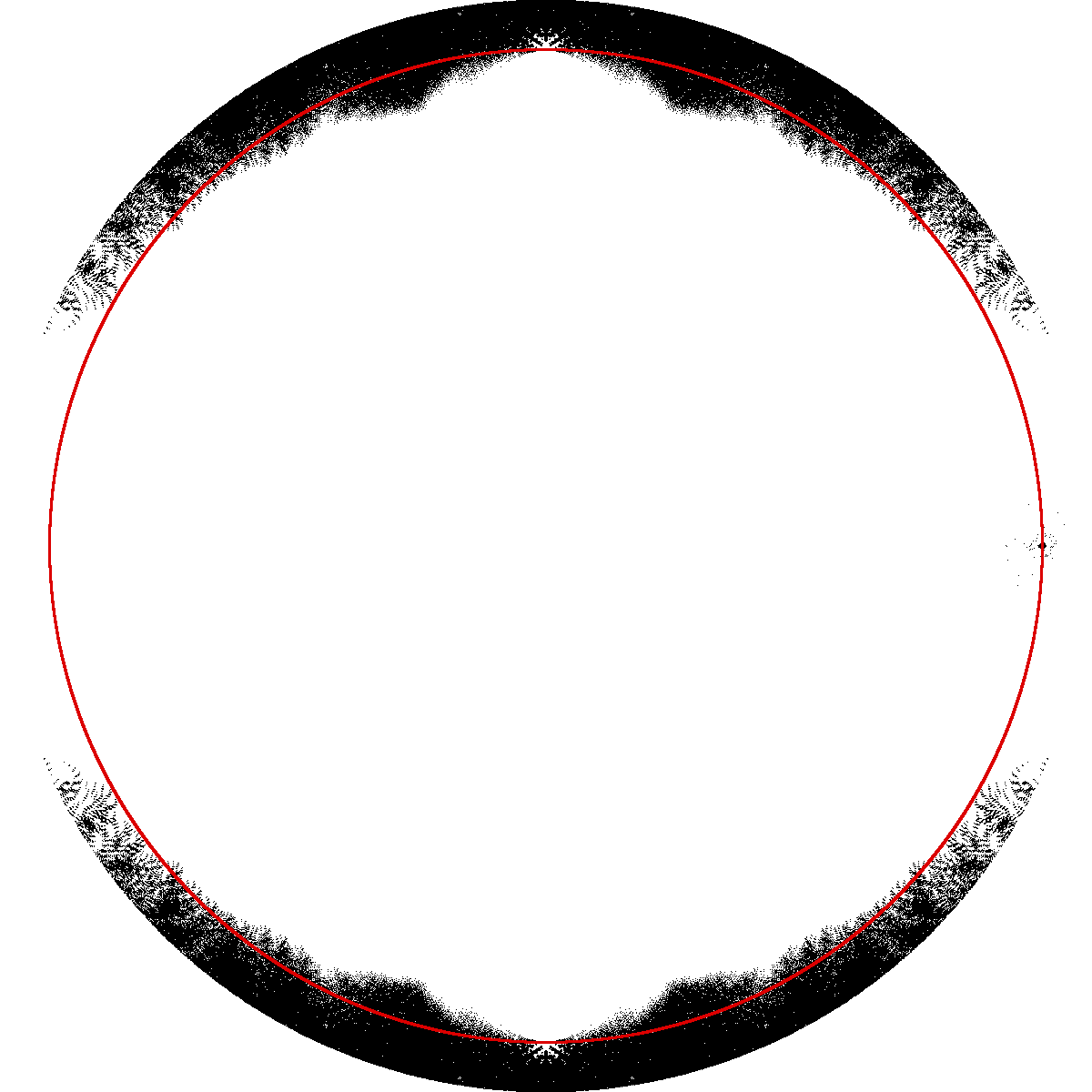}
\end{minipage}
\hfill
\begin{minipage}[c]{0.48\textwidth}
\centering
\includegraphics[width=0.909\textwidth]{W.png}
\end{minipage}
\caption{On the left, the roots of the differences $p_u-p_v$ for all words $u,v$ of
length at most $12$ which lie in $|z-1/2|<0.55$; the circle $|z-1/2|=1/2$ is drawn
in red. On the right, $\W$ (in white) as in Figure~\ref{W}, at the same scale.
The roots accumulate on $\RR$.}\label{roots}
\end{figure}

Since $z$ lies in $\gamma$, the point $p_u(z)=u(z)$ lies in $u\gamma$. Hence if
$u\in fS$ and $v\in gS$ satisfy $p_u(z)=p_v(z)\ne z$ then $f\gamma\cap g\gamma$
contains a point other than $z$, and $z\in \RR$ by Lemma~\ref{lemma:z_in_W}. So
every root of a difference $p_u-p_v$ which lies in the disk $|z-1/2|<1/2$, and at
which the common value is not $z$, belongs to $\RR$; and $\RR$ is the closure of the
set of all such roots. For the reverse inclusion one truncates at length $n$ the two
addresses of a point $p\ne z_0$ of $f\gamma\cap g\gamma$, obtaining $u_n\in fS$ and
$v_n\in gS$ with $|p_{u_n}(z_0)-p_{v_n}(z_0)|$ of order
$\max(|z_0|,|1-z_0|)^n$, and replaces these near-roots by genuine roots converging
to $z_0$; we omit the details, which we do not need.
Figure~\ref{roots} shows these roots for all $u,v$ of length at most $12$, beside
Figure~\ref{W} at the same scale.

The parallel with $\mathcal{M}$ is close. There $\lambda \in \mathcal{M}$ if and only
if $\sum_n a_n\lambda^n = 0$ for some sequence $a_n \in \lbrace 0,\pm 1\rbrace$, and
$\mathcal{M}$ is the closure of the roots in the unit disk of the corresponding
polynomials; the pictures of Odlyzko and Poonen \cite{Odlyzko_Poonen} are of exactly
this kind. Here the $\lbrace 0,\pm 1\rbrace$ polynomials are replaced by the
differences $p_u-p_v$.

The roots of $p_u-p_v$ are not confined to the disk $|z-1/2|<1/2$, and only those
inside it are relevant. The excluded pair of \S~\ref{subsection:certifying_W}
illustrates this:
$$p_{ffg}(z)-p_{ggf}(z) = z(z-1)(2z^2-2z+1),$$
whose roots $0$, $1$ and $(1\pm i)/2$ all lie on the circle $|z-1/2|=1/2$.

\subsection{Certifying $z$ in the interior of $\RR$}

We now show how to modify the algorithm from \S~\ref{subsection:certifying_W}
to certify (numerically) that a point $z$ is in the interior of $\RR$. One consequence of the
nature of the algorithm is that it implies that the interior of $\RR$ is
dense in $\RR$.

The idea is a modification of the method of {\em traps}, introduced in
\cite{Calegari_Koch_Walker} to prove Bandt's Conjecture on interior points
in $\mathcal{M}$, the connectivity locus for another 1-parameter
family of complex 1-dimensional IFSs (and the analog to $\RR$ for this family)

The idea is very simple. Let $u\in fS$ and $v\in gS$ and choose some finite $n$.
Let $uS_n B$ denote the union of translates of $B$ for all elements of $uS_n$
and likewise $vS_n B$.

We will shortly give the definition of a {\em stable crossing}, which certifies an
interior point in $\RR$, but first we clarify two points of terminology. A {\em ray}
from a point $p$ means a proper embedding $r:[0,\infty) \to \C$ with $r(0)=p$; it is a
topological ray, and need not be a Euclidean half-line. And if $r^\pm$, $s^\pm$ are
four pairwise disjoint rays, we say that $r^\pm$ {\em link} $s^\pm$ at infinity if the
following holds. Choose a closed round disk $D$ containing $uS_nB \cup vS_nB$ in its
interior, and truncate each of the four rays at the first time it meets $\partial D$.
Since the rays are disjoint this produces four distinct points of $\partial D$,
together with four paths in $D$ from the initial points of the rays to them; we
require that for some such $D$ these four points occur in the cyclic order
$r^+,s^+,r^-,s^-$ on $\partial D$.

\begin{definition}[Stable Crossing]\label{def:crossing}
The triple $(u,v,n)$ is a {\em stable crossing} if there
are {\em disjoint proper rays} $r^\pm$, $s^\pm$ from points $p^\pm$ and $q^\pm$
in $u\gamma$, $v\gamma$
to infinity so that $r^\pm$ are disjoint from $vS_nB$, so that $s^\pm$ are
disjoint from $uS_nB$, and so that $r^\pm$ link $s^\pm$ at infinity.
\end{definition}

The existence of a stable crossing is, in fact, stable in $z$; and given
$(u,v,n)$ stable for $z$ one may estimate a lower bound on the radius of
a ball around $z$ for which this triple continues to be stable.
We obtain an explicit bound of this kind in \S~\ref{subsection:proof}, in the
course of describing the computer implementation.

Furthermore, the existence of a stable crossing certifies $z\in \RR$ for
elementary topological reasons.

The topological input is the following standard fact.

\begin{lemma}[Linking]\label{lemma:linking}
Let $D$ be a closed disk and let $a^+,b^+,a^-,b^-$ be four distinct points of
$\partial D$ occurring in this cyclic order. Let $A$ be a path in $D$ from $a^+$ to
$a^-$, and let $B$ be a path in $D$ from $b^+$ to $b^-$. Then
$A\cap B \ne \emptyset$.
\end{lemma}

Neither $A$ nor $B$ is required to be embedded. After a homeomorphism of $D$ to a
square carrying the four points to the midpoints of the four sides, this is the
familiar statement that a path joining the left side to the right side must meet a
path joining the top to the bottom; it is equivalent to the Jordan curve theorem, and
follows from Brouwer's theorem applied to $(s,t) \to A(s)-B(t)$.

\begin{lemma}
If $(u,v,n)$ is a stable crossing for $z$ then $z\in \RR$.
\end{lemma}
\begin{proof}
Let $p^\pm \in u\gamma$ and $q^\pm \in
v\gamma$ be the finite endpoints of $r^\pm$ and $s^\pm$ respectively.
Since $u\gamma$ and $v\gamma$ are path connected, there are arcs
$\alpha \subset u\gamma$ and $\beta \subset v\gamma$ joining $p^\pm$ and
$q^\pm$. Then $r^+\cup \alpha \cup r^-$ and $s^+\cup \beta\cup s^-$ are
properly immersed lines in $\C$ that are embedded and disjoint at infinity
where their endpoints are linked. Thus their algebraic
intersection number (rel. their ends) is odd, so they must intersect.
Concretely: truncating the four rays at $\partial D$ as above and applying
Lemma~\ref{lemma:linking} to $A = r^+\cup \alpha\cup r^-$ and
$B = s^+\cup \beta \cup s^-$ shows $A\cap B \ne \emptyset$.
But by hypothesis the only place they might intersect is $\alpha$ with
$\beta$.

It remains to check that a point of $\alpha \cap \beta$ is not the point $z$, which
lies in $f\gamma\cap g\gamma$ for every $z$ and so certifies nothing. Suppose then
that $w_z$ is injective, and for a word $w$ let $I_w\subset I$ denote the interval of
parameters whose address begins with $w$, so that $w_z(I_w) = w\gamma$. Since $u$
begins with $f$ and $v$ begins with $g$ we have $I_u \subset [0,1/2]$ and
$I_v \subset [1/2,1]$, so that $I_u \cap I_v \subset \lbrace 1/2\rbrace$; as $w_z$ is
injective it follows that $u\gamma\cap v\gamma$ is either empty or equal to
$\lbrace z \rbrace$. It is not empty, since it contains $\alpha\cap\beta$; so
$\alpha\cap\beta = \lbrace z\rbrace$ and $1/2$ is the right hand endpoint of $I_u$.
Then $z = w_z(1/2)$ is an endpoint of the arc $u\gamma$, and $\alpha$ is a subarc of
$u\gamma$ containing $z$, so $z$ is an endpoint of $\alpha$; that is,
$z \in \lbrace p^+,p^-\rbrace$. But $p^\pm$ lie on $r^\pm$, which are disjoint from
$vS_nB \supset v\gamma$, whereas $z \in v\gamma$ --- a contradiction. So $w_z$ is not
injective, and $z\in \RR$ by Lemma~\ref{lemma:z_in_W}.
\end{proof}

To search for stable crossings we enumerate pairs $(u,v)$ by the algorithm
from the previous section, then for each we compute $uS_nB$ and $vS_nB$ for
some finite $n$ and look for rays $r^\pm$, $s^\pm$ and points $p^\pm$, $q^\pm$ 
forming a stable crossing. If we find one we certify $z$ (and some ball of
computable radius about it) as lying in $\RR$.

\medskip

Let's say more generally that two compact path-connected subsets 
$K,L \subset \C$ have a stable crossing if there is some $\epsilon>0$
and disjoint rays $r^\pm$, $s^\pm$
from points $p^\pm$, $q^\pm$ in $K$ and $L$ so that $r^\pm$ is disjoint
from the $\epsilon$-neighborhood of $L$ and $s^\pm$ is disjoint from the
$\epsilon$-neighborhood of $K$, and $r^\pm$, $s^\pm$ link at infinity.
Thus $u\gamma,v\gamma$ have a stable crossing if and only if
$(u,v,n)$ is a stable crossing for some $n$. See Figure~\ref{Gamma_crossings}
for some relevant examples.

As an application of stable crossings we prove that the interior of $\RR$ is
dense in $\RR$. Both the statement and the proof are very similar to 
Theorem~7.2.7 from \cite{Calegari_Koch_Walker}.

\begin{theorem}[Interior is Dense]\label{interior_is_dense}
The interior of $\RR$ is dense in $\RR$.
\end{theorem}
\begin{proof} 
It suffices to find a stable crossing arbitrarily close to any $z_0\in \RR$.
Let $z_0\in \RR$ so that $u\gamma$ intersects $v\gamma$ for some finite
words $u,v$ which do not start with the prefix $ffg,ggf$. Without loss of
generality, we may assume $u,v$ both have length $n\gg 1$. There are
complex numbers $\alpha,\beta$ so that $\gamma = \alpha v^{-1}u\gamma+\beta$, and
by padding $u$ or $v$ or both with additional letters $f$ or $g$ if necessary we
may arrange for $1/C<|\alpha|<C$ for some $C$ depending on $z$ while still having
$u\gamma$ intersect $v\gamma$. 

Actually, with more work we may eliminate the factor of $\alpha$.
Note that $\alpha(z)$ is just a product of integer powers of $-z$ and $(z-1)$
according to how many $f$s and $g$s are in $u$ and $v$. Thus by multiplying
$u$ and $v$ on the right by suitable powers of $f$ and $g$ we may arrange for
$\gamma(z) = v^{-1}u\gamma(z) + \beta(z)$ for some new $\beta(z)$.
It follows that to prove the theorem we just need to show that
{\em whatever} shape $\gamma(z_0)$ is, there is some constant $\mu\in \C$ so that
$\gamma(z_0)$ and $\gamma(z_0)+\mu$ have a stable crossing. This is proved in
Lemma~\ref{translation_lemma} and Lemma~\ref{gamma_not_convex}.

We now give this argument in detail, following the proof of Theorem~7.2.7 of
\cite{Calegari_Koch_Walker}. For an infinite word $a\in \lbrace f,g\rbrace^\infty$
let $P(a,z)$ denote the point of $\gamma(z)$ with address $a$, that is the limit of
$a_1a_2\cdots a_m(x)$ as $m\to \infty$; the limit exists and is independent of $x$
because $f$ and $g$ are contractions, and since each $a_1\cdots a_m(x)$ is a
polynomial in $z$ and the convergence is locally uniform, $P(a,z)$ is holomorphic in
$z$. Since $z_0\in \RR$ there are infinite words $a$ beginning with $f$ and $b$
beginning with $g$ with $P(a,z_0)=P(b,z_0)$, and by the discussion in
\S~\ref{subsection:certifying_W} we may take $(a,b)$ to be other than the pair
$(f(fg)^\infty, g(gf)^\infty)$ of addresses of the parameter $1/2$. Put
$$T(z) := P(b,z)-P(a,z),$$
a holomorphic function on the disk with $T(z_0)=0$.

The function $T$ is not identically zero. For $\gamma(1/2)=[0,1]$ and $w_{1/2}$ is
injective, so $T(1/2)=0$ would force $a$ and $b$ to be addresses of the same
parameter; since $a$ begins with $f$ and $b$ with $g$ that parameter would be $1/2$,
whose two addresses are $f(fg)^\infty$ and $g(gf)^\infty$ --- the pair we excluded.
So $T$ is a nonconstant holomorphic function vanishing at $z_0$, and by the open
mapping theorem $T(U)$ contains the disk $|\mu|<2\delta$ for some $\delta>0$, where
$U$ is the disk of radius $\epsilon$ about $z_0$.

Let $u_m,v_m$ be the prefixes of $a,b$ of length $m$, and let $d_m$ be the number of
$f$s in $u_m$ minus the number in $v_m$. Set $U_m := u_mg^{d_m}$ and
$V_m := v_mf^{d_m}$ if $d_m\ge 0$, and $U_m := u_mf^{-d_m}$, $V_m := v_mg^{-d_m}$
otherwise. Then $U_m$ and $V_m$ begin with $f$ and $g$ and contain equally many $f$s
and equally many $g$s, so their linear parts coincide; call the common value
$\alpha_m$, and note $|\alpha_m| \to 0$. Since $p_{U_m}(z)$ and $P(a,z)$ both lie in
$u_m\gamma(z)$, whose diameter is $|\alpha_{u_m}|\textnormal{diam}(\gamma(z))$, we get
$$|T_m(z)-T(z)| \le \left( |\alpha_{u_m}|+|\alpha_{v_m}| \right)
\textnormal{diam}(\gamma(z)) \longrightarrow 0$$
uniformly on $U$, where $T_m(z) := p_{V_m}(z)-p_{U_m}(z)$. By Rouch\'e it follows
that $T_m(U)$ contains the disk $|\mu|<\delta$ for all sufficiently large $m$.
Finally $U_m^{-1}V_m$ is the translation by $\beta_m := T_m/\alpha_m$, so
$U_m\gamma(z)$ meets $V_m\gamma(z)$ exactly when $\gamma(z)$ meets
$\gamma(z)+\beta_m(z)$; and $\beta_m(U)$ contains the disk of radius
$\delta/|\alpha_m|$ about $0$, which for large $m$ contains any prescribed $\mu$.

Two points remain. Lemma~\ref{translation_lemma} applies to compact {\em full} sets,
and $\gamma(z_0)$ is not full when $w_{z_0}$ is not injective; it should be applied
instead to the {\em full set} $X(z)$, the union of $\gamma(z)$ with the bounded
components of its complement, which is the smallest compact full set containing
$\gamma(z)$. This is harmless, because the four points witnessing the crossing of
$X(z_0)$ and $X(z_0)+\mu$ lie on the outermost boundary of
$X(z_0)\cup (X(z_0)+\mu)$, and $\partial X(z_0) \subset \gamma(z_0)$; so the arcs and
rays may be taken in $\gamma(z_0)$ and $\gamma(z_0)+\mu$, and these have a stable
crossing. Compare Lemma~7.2.3 of \cite{Calegari_Koch_Walker}, which characterizes
those $z$ for which the filling of the limit set is convex.
\end{proof}

Thus Theorem~\ref{interior_is_dense} is reduced to the following two lemmas:

\begin{lemma}[Translation Lemma]\label{translation_lemma}
Let $K$ be a compact full subset of $\C$ and suppose $K$ is not convex. 
There is some complex number $\mu$ of order $\text{diam}(K)$
so that $K$ and $K + \mu$ have a stable crossing.
\end{lemma}

This is Lemma~7.2.2 from \cite{Calegari_Koch_Walker}. Finally we need to understand
the set of $z$ for which $X(z)$ is convex.

\begin{lemma}[$X$ not convex]\label{gamma_not_convex}
For $|z-1/2|<1/2$ and $z$ not real, $X(z)$ is not convex.
\end{lemma}
Of course if $|z-1/2|<1/2$ and $z$ {\em is} real then $\gamma(z)=[0,1]$ and
$z\in \W$.
\begin{proof}
Write $X:=X(z)$ and suppose $X$ is convex, with $z$ not real. Since
$z=w_z(1/2)$ lies in $\gamma(z)\subset X$ and is not real, $X$ is not contained in a
line; being convex it therefore has
nonempty interior (in particular $X$ has Hausdorff dimension $2$). We claim
that actually $\gamma(z)=X$. To see this, let's let $p\in X$ be a
point which realizes the maximal distance $\epsilon$ to $\gamma(z)$. Because $X$
is convex, and $f(X)$, $g(X)$ intersect, it follows that $X=f(X)\cup g(X)$
(without convexity there might be some omitted region `trapped' between $f(X)$ and
$g(X)$). But then $p\in f(X)$ or $g(X)$; thus the distance from $p$ to
$\gamma(z)$ is at most $\max(|z|\epsilon,|1-z|\epsilon)$. Since $|z|$ and
$|1-z|$ are both strictly less than $1$ it follows that $\epsilon=0$ and
therefore $\gamma(z)=X$ and has Hausdorff dimension $d=2$. But then
$|z|^2 +|1-z|^2\ge 1$ so $|z-1/2|\ge 1/2$ and we are done.
\end{proof}

Compare to Lemma~7.2.3 from \cite{Calegari_Koch_Walker}.
This completes the proof of Theorem~\ref{interior_is_dense}.

\subsection{Proof of Theorem~\ref{theorem:wiggle_island}}\label{subsection:proof}

Using the techniques of the previous two sections, we may numerically certify
some $z\in \W$ and numerically check that a polygonal loop separating $z$ from
$0$ is contained in $\RR$. To make this computationally feasible we pursue
the following strategy:
\begin{enumerate}
\item{find an apparent island by some experimentation, in a small region
$U\subset \C$;}
\item{fix a polygon $\Gamma$ and verify that $\gamma(z)\subset \Gamma$ for all
$z\in U$;}
\item{compute a region $O\subset \C^*\times \C$ so that $\gamma(z)$ and
$\alpha \gamma(z)+\beta$ have a stable crossing for all $(\alpha,\beta)\in O$
and all $z\in U$;}
\item{run the algorithm at parameters $z\in U$ to generate crossing
pairs $(u,v)$ and for each compute $(\alpha(z),\beta(z))$ such that
$\gamma(z) = \alpha(z) v^{-1}u\gamma(z) + \beta(z)$;}
\item{if $(\alpha(z),\beta(z))$
is in $O$, compute an $\epsilon$ so that $(\alpha(z'),\beta(z'))\in O$ when
$|z'-z|<\epsilon$; and finally}
\item{cover a circle surrounding the island by such $\epsilon$-disks, and
certify that its centre lies in $\W$; this proves the theorem.}
\end{enumerate}

We sketch these points here; the computation itself, and the data it produces, are
described in \S~\ref{subsection:computation}.

\begin{figure}[hpbt]
\centering
\includegraphics[scale=0.2]{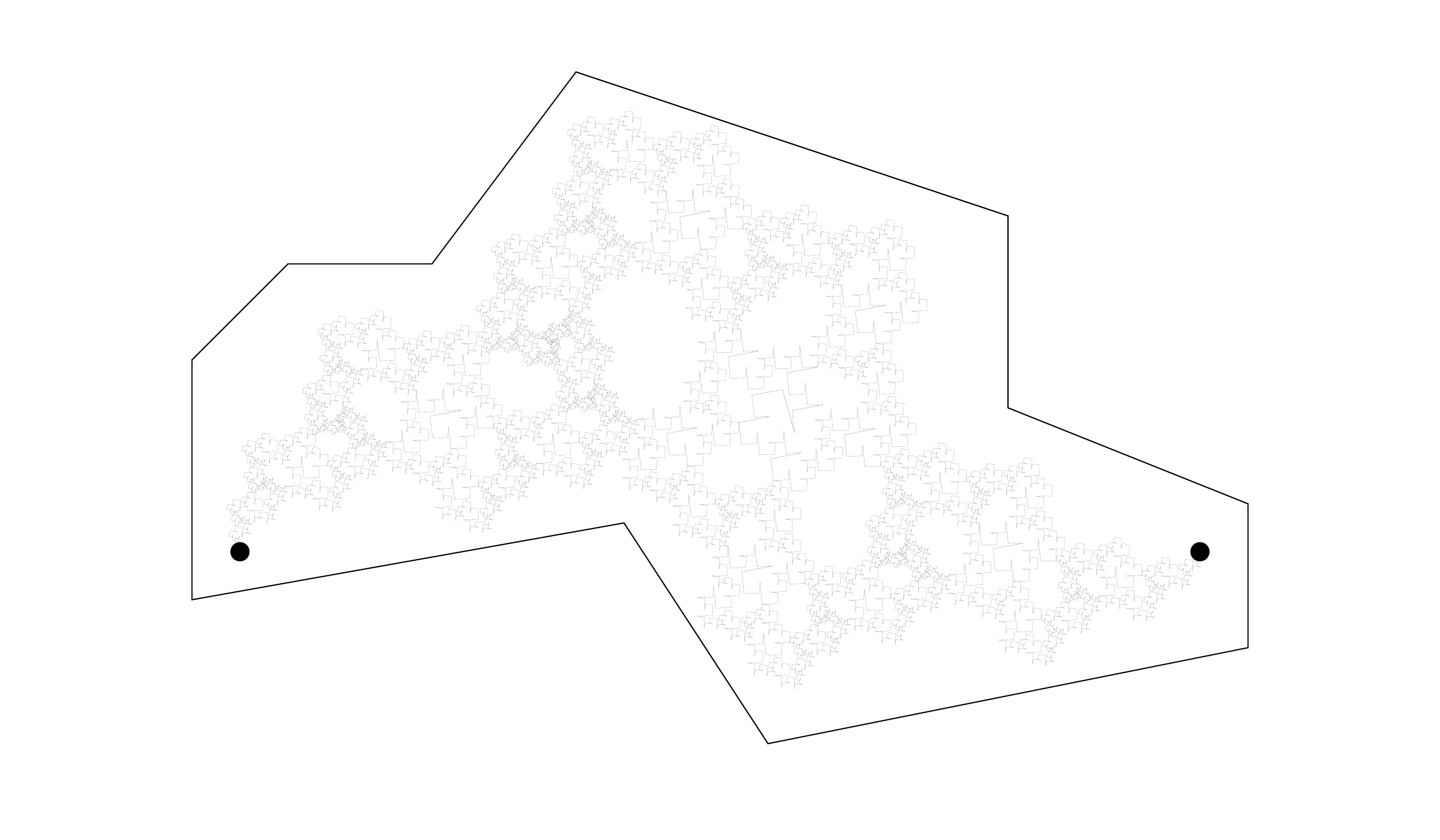}
\caption{The disk $\Gamma$ encloses $\gamma(z)$ for all $z\in U$. The vertices
$0$ and $1$ of $\gamma(z)$ are distinguished.}\label{Gamma}
\end{figure}

The region $U$ is the disk $|z-(0.3409+0.43486i)|\le 10^{-4}$.
We may then readily determine the vertices of a polyhedral disk $\Gamma$
that contains a thin neighborhood of $\gamma(z)$ for all $z\in U$; see
Figure~\ref{Gamma}. The containment is not merely apparent: it is verified in
\S~\ref{subsection:computation}.

\begin{figure}[hpbt]
\centering
\includegraphics[scale=0.7]{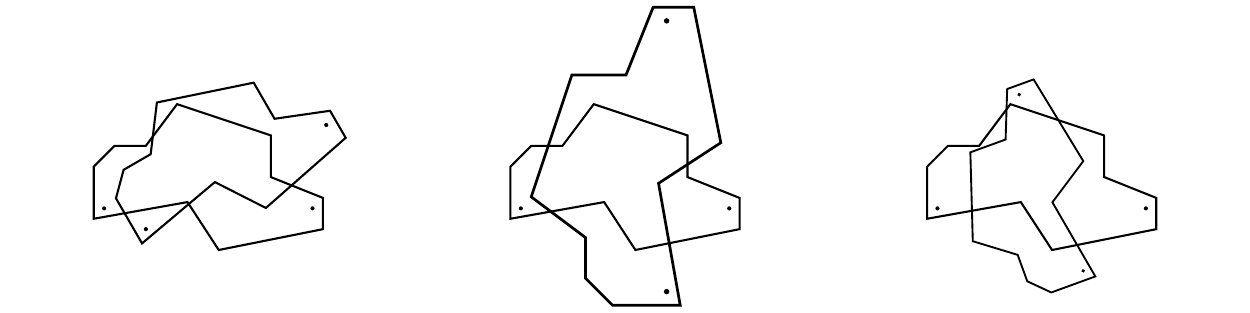}
\caption{Stable crossings of $\Gamma$ and $\alpha\Gamma+\beta$.}\label{Gamma_crossings}
\end{figure}

Next we compute some open polydisks in $\C^*\times \C$ parameterizing 
$(\alpha,\beta)$ for which $\Gamma$ and $\alpha \Gamma + \beta$ 
have a suitable stable crossing (with rays landing at $0,1$ and $\beta,\alpha+\beta$
respectively, which are always in $\gamma(z)$ and $\alpha\gamma(z)+\beta$);
see Figure~\ref{Gamma_crossings}.

\begin{figure}[hpbt]
\centering
\includegraphics[width=0.5\textwidth]{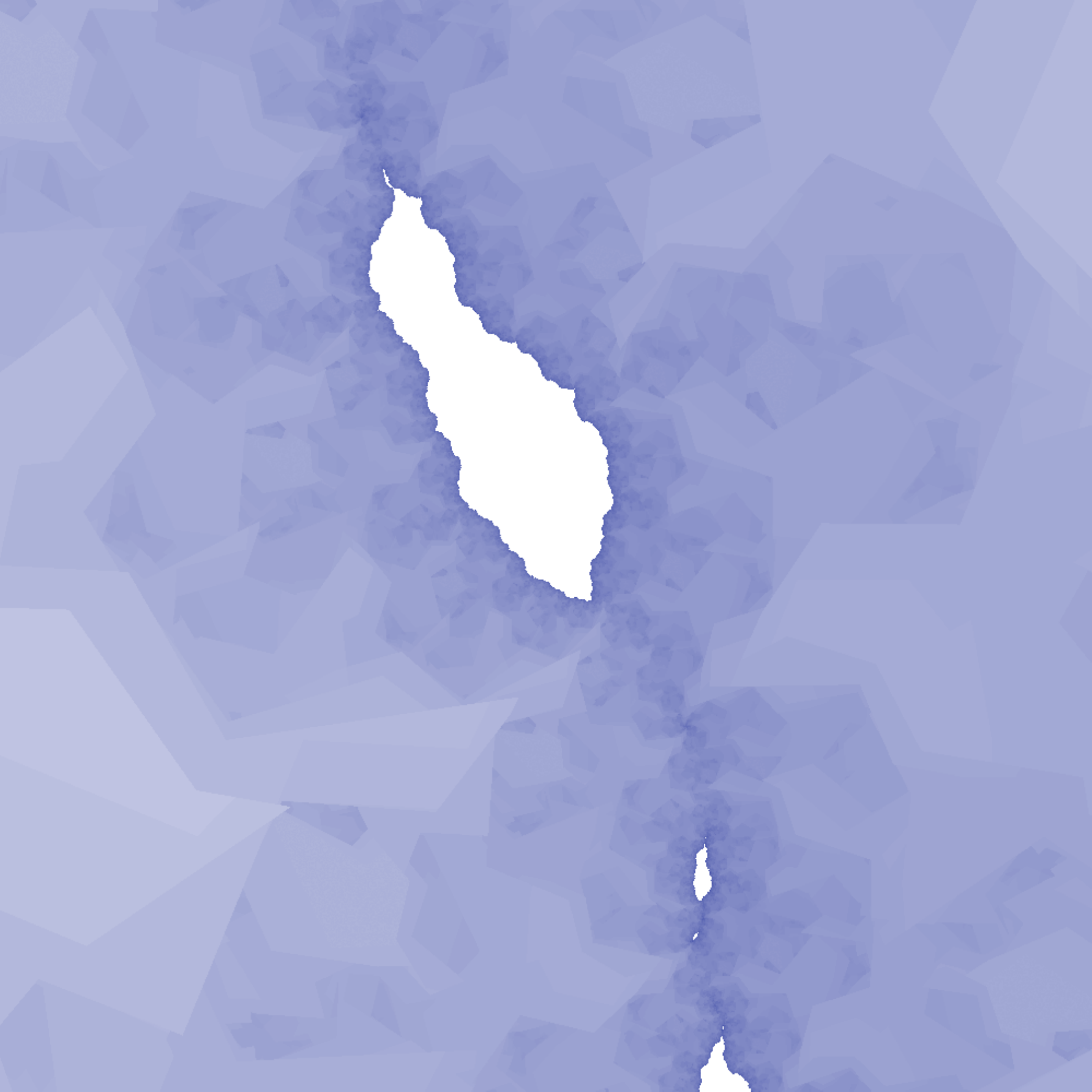}
\caption{Certificates on the square of half width $5\times 10^{-5}$ about
$0.3409+0.43486i$: white is certified to lie in $\W$, blue in $\RR$, as in
Figure~\ref{W}. No parameter of the grid is left undecided.}\label{W_island_rigorous}
\end{figure}

The program {\tt wiggle} implements this strategy. Figure~\ref{W_island_rigorous}
shows the result of running it at every point of a grid near the island: at each
parameter the algorithm terminates, either with a stable crossing certifying
$z\in \RR$, or with an empty list of viable pairs certifying $z \in \W$. Wiggle
Island is the elongated white region, with a smaller white region to the southeast.

Intermediate between Wiggle Island and the mainland there is apparent in the Figure
a smaller island, and a speck (which on magnification turns out to be an island
too).
It seems likely that there are infinitely many islands,
arranged in an asymptotically geometric spiral converging to an algebraic point
$z_c\sim 0.340922 + 0.43481i$, though we have not attempted to prove this. The
analogous phenomenon for $\mathcal{M}$ --- a spiral of islands accumulating on a
point of the frontier --- is certified in \S~9 of \cite{Calegari_Koch_Walker}, and we
refer the interested reader to that paper.

We end with a conjecture, parallel to Conjecture~9.2.7 from \cite{Calegari_Koch_Walker}:

\begin{conjecture}
Every point in the frontier of $\RR$ is the limit of islands in $\W$ of diameter
going to zero.
\end{conjecture}

\subsection{The computation}\label{subsection:computation}

We describe the computation in enough detail that it may be checked. The program
{\tt wiggle} is available with this paper, and the certificate data discussed below
is in the file {\tt certified\_loop.txt}.

\medskip\noindent
{\em The envelope.} Let $\Gamma$ be the closed polygon with vertices
$$\begin{array}{llll}
(0.35,0.5), & (0.2,0.3), & (0.05,0.3), & (-0.05,0.2),\\
(-0.05,-0.05), & (0.4,0.03), & (0.55,-0.2), & (1.05,-0.1),\\
(1.05,0.05), & (0.8,0.15), & (0.8,0.35). &
\end{array}$$
Every crossing certificate below depends on the containment
$\gamma(z)\subset \Gamma$, which we verify directly. Since
$\gamma = \bigcup_{|w|=k} w\gamma$ and $\gamma\subset B$, we have
$$\gamma(z) \subset \bigcup_{|w|=k} wB,$$
and $wB$ is the disk of radius $|\alpha_w|R$ about $w(1/2)$, where $\alpha_w$ is the
linear part of $w$. Evaluating $w(1/2)$ and $\alpha_w$ in disk arithmetic over a disk
of parameters, and bounding $R$ above there, reduces the containment to finitely many
inequalities. For the parameter disk $|z-(0.3409+0.43486i)|\le 10^{-4}$ one gets
$R\le 1.3142$, and at $k=18$ all $2^{18}=262144$ disks $wB$ lie in $\Gamma$, the
smallest clearance from $\partial\Gamma$ being $0.0042$. (At $k=10$ and $k=14$ the
covering disks are still too large and the check fails.) Every parameter used below
lies in this disk.

\medskip\noindent
{\em The crossing test.} For a pair of words $u\in fS$, $v\in gS$ let
$A=u^{-1}v: w\to \alpha w+\beta$, so that $u\gamma$ meets $v\gamma$ exactly when
$\gamma$ meets $A\gamma$. To certify a stable crossing we exhibit the four rays
themselves. Their feet are the four endpoints
$$0,\ 1 \in \gamma, \qquad A(0)=\beta,\ A(1)=\alpha+\beta \in A\gamma,$$
and we ask for four {\em straight} rays from these points whose arguments occur in
alternating cyclic order, such that the rays from $0$ and $1$ are disjoint from
$A\Gamma$, the rays from $A(0)$ and $A(1)$ are disjoint from $\Gamma$, and each ray
of the first pair is disjoint from each ray of the second. When such rays exist the
two proper paths $r^+\cup\gamma\cup r^-$ and $s^+\cup A\gamma\cup s^-$ can meet only
along $\gamma\cap A\gamma$, and Lemma~\ref{lemma:linking} forces that intersection to
be nonempty. In practice we first try the four rays emanating radially from the point
where the segments $[0,1]$ and $[A(0),A(1)]$ cross; these automatically have
alternating arguments and are pairwise disjoint, so only disjointness from the
envelopes need be checked. If that fails we search a grid of directions.

\medskip\noindent
{\em Certified radii.} Fix the four directions at the values they take at a parameter
$z$. Then the cyclic order cannot change as $z$ varies, and the certificate persists
as long as eight clearances remain positive: the distances from the rays at $0$ and
$1$ to $A\Gamma$, the distances from the rays at $A(0)$ and $A(1)$ to $\Gamma$, and
the four distances between rays of opposite type. Each clearance is $1$--Lipschitz in
the objects that move: $A\Gamma$ moves at rate at most
$\max_v |v| \cdot |\alpha'| + |\beta'|$, where $v$ runs over the vertices of $\Gamma$,
the foot $A(0)$ at rate $|\beta'|$, and the foot $A(1)$ at rate
$|\alpha'|+|\beta'|$. Writing $m_i$ for the clearances and $L_i$ for these rates, the
crossing is certified on the disk of radius
$$\epsilon = \min_i m_i/L_i.$$
We bound $|\alpha'|$ and $|\beta'|$ above on that disk by carrying $d\alpha/dz$ and
$d\beta/dz$ along the composition defining $A$ and evaluating in disk arithmetic, so
that they are genuine upper bounds and not merely derivatives at the centre; the disk
is enlarged until it contains the resulting $\epsilon$. Because $|\alpha|^{-1}$ is
exponentially large in the length of $u$ and $v$, short words give much larger
$\epsilon$, and the program certifies every crossing available over several levels of
the search and keeps the largest radius.

\begin{figure}[hpbt]
\centering
\begin{tikzpicture}[x=0.62\textwidth,y=0.62\textwidth,
	loopdisk/.style={fill=orange!85!red,draw=none},
	loopcircle/.style={draw=black,line width=0.3pt},
	loopcentre/.style={fill=red!65!black,draw=none}]
\node[anchor=south west,inner sep=0] at (0,0)
	{\includegraphics[width=0.62\textwidth]{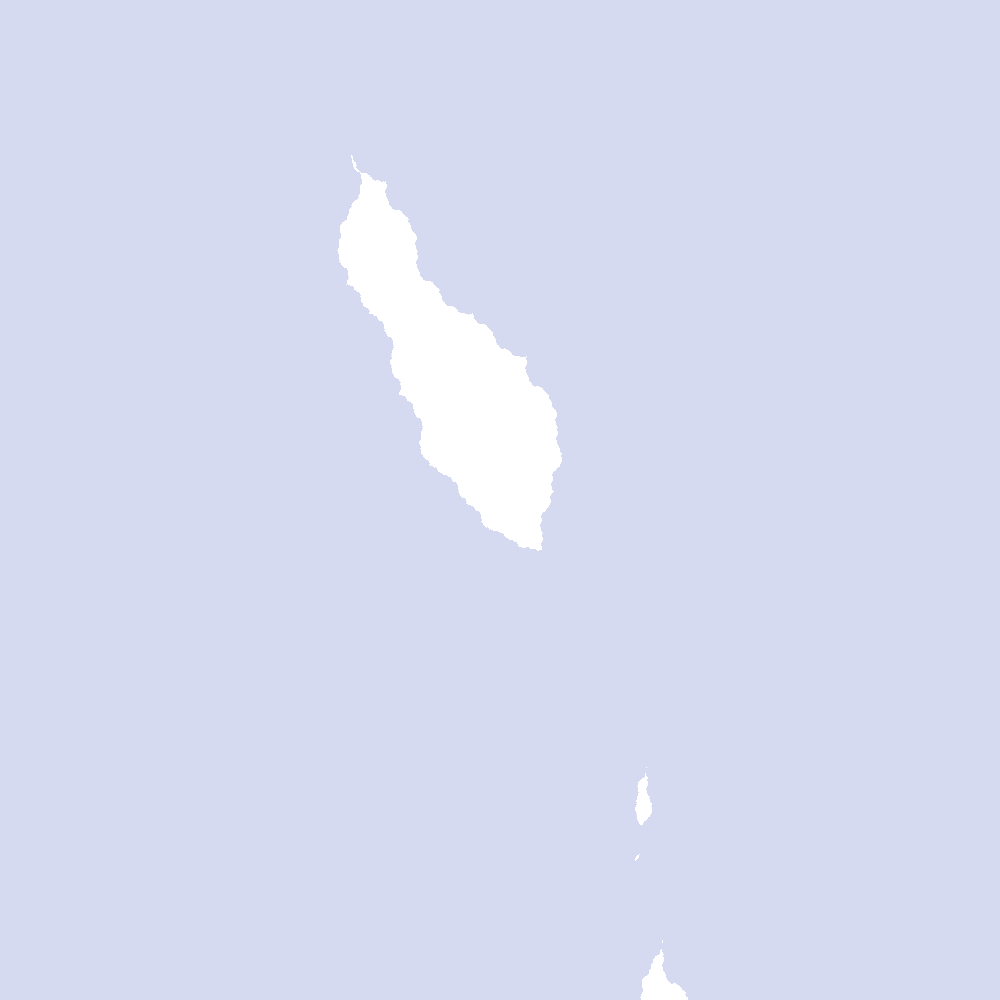}};
\input{certified_loop_tikz.tex}
\end{tikzpicture}
\caption{The certified loop. Pale blue is certified to lie in $\RR$ and white in
$\W$; the band of orange disks are the $1046$ disks, each certified to lie in
$\RR$, which cover the circle $|z-L|=4.2\times 10^{-5}$ drawn in black; and the dark
red dot at its centre is $L$ itself, certified to lie in $\W$. The window has half
width $5\times 10^{-5}$.}\label{certified_loop}
\end{figure}

\medskip\noindent
{\em The loop.} Let $L:=0.34090875+0.43484625i$. The circle $|z-L|=4.2\times 10^{-5}$
is covered by $1046$ certified disks, with no gaps; the certified radii range from
$3.5\times 10^{-10}$ to $3.2\times 10^{-6}$, with median $2.0\times 10^{-7}$, and only
$202$ distinct pairs $(u,v)$ occur, of lengths between $21$ and $43$. Rather than
walking the circle, the program samples a coarse set of angles --- each certified
disk covers an arc $|\theta-\theta_i|<2\arcsin(\epsilon_i/2\,\textnormal{rad})$ ---
and then bisects only the angular gaps that remain, so that effort is spent where it
is needed. See Figure~\ref{certified_loop}.

Finally the algorithm of \S~\ref{subsection:certifying_W} terminates at $z=L$ itself,
certifying $L\in \W$. Since the circle lies in $\RR$, the component of $\W$ containing
$L$ is contained in the disk bounded by that circle, and so is not the component
containing $1/2$. This proves Theorem~\ref{theorem:wiggle_island}.

\medskip\noindent
{\em Caveats.} The disk arithmetic inflates radii by a fixed relative factor to
absorb rounding, which is not the same as directed rounding; and the certification of
$L\in \W$, which rests on the disjointness test for pairs of balls, is carried out in
ordinary floating point. Neither is likely to be delicate --- the margins above are
many orders of magnitude larger than double precision --- but a fully rigorous
implementation would use interval arithmetic throughout.

\begin{remark}
In fact the existence of Wiggle Island has been well-known to Europeans at least
since the early 17th century, and is well-documented in the literature \cite{Wiggles}.
\end{remark}


\begin{thebibliography}{99}
\bibitem{Aseev_Tetenov_Kravchenko}
	V. Aseev, A. Tetenov and A. Kravchenko,
	\emph{On self-similar Jordan curves on the plane},
	Siberian Math. J. {\bf 44} (2003), no. 3, 379--386
\bibitem{Astala}
	K. Astala,
	\emph{Area distortion of quasiconformal mappings},
	Acta Math. {\bf 173} (1994), no. 1, 37--60
\bibitem{Astala_Rohde_Schramm}
	K. Astala, S. Rohde and O. Schramm, unpublished work; see the lecture of
	S. Rohde, \emph{Jordan curves and dimension of quasicircles}, at the
	Oded Schramm Memorial Conference,
	\url{https://www.microsoft.com/en-us/research/video/oded-schramm-memorial-conference-day-two-session-two/}
\bibitem{Bandt}
	C. Bandt,
	\emph{On the Mandelbrot set for pairs of linear maps},
	Nonlinearity {\bf 15} (2002), no. 4, 1127--1147
\bibitem{Barnsley_Harrington}
	M. Barnsley and A. Harrington,
	\emph{A Mandelbrot set for pairs of linear maps},
	Phys. D {\bf 15} (1985), no. 3, 421--432
\bibitem{Bousch}
	T. Bousch,
	\emph{Paires de similitudes $Z \to sZ+1$, $Z \to sZ-1$},
	unpublished manuscript (1988),
	\url{https://www.imo.universite-paris-saclay.fr/~thierry.bousch/preprints/paires_sim.pdf}
\bibitem{Calegari_Koch_Walker}
	D. Calegari, S. Koch and A. Walker,
	\emph{Roots, Schottky semigroups, and a proof of Bandt's Conjecture},
	Ergodic Theory and Dynamical Systems {\bf 37} (2017), no. 8, 2487--2555
\bibitem{Connelly_Demaine_Rote}
	R. Connelly, E. Demaine and G. Rote,
	\emph{Straightening polygonal arcs and convexifying polygonal cycles},
	Proc. 41st Annual Symposium on Foundations of Computer Science, 432--442,
	Redondo Beach, California, November 2000.
\bibitem{de_Rham_1956}
	G. de Rham,
	\emph{Sur une courbe plane},
	J. Math. Pures Appl. (9) {\bf 35} (1956), 25--42
\bibitem{de_Rham_1957}
	G. de Rham,
	\emph{Sur quelques courbes d\'efinies par des \'equations fonctionnelles},
	Rend. Sem. Mat. Univ. Politec. Torino {\bf 16} (1957), 101--113
\bibitem{Odlyzko_Poonen}
	A. Odlyzko and B. Poonen,
	\emph{Zeros of polynomials with $0,1$ coefficients},
	Enseign. Math. (2) {\bf 39} (1993), no. 3-4, 317--348
\bibitem{Solomyak}
	B. Solomyak,
	\emph{On the `Mandelbrot set' for pairs of linear maps: asymptotic
	self-similarity}, Nonlinearity {\bf 18} (2005), no. 5, 1927--1943
\bibitem{Wiggles}
	Wikipedia page \url{https://en.wikipedia.org/wiki/The_Wiggles}
\end{thebibliography}
\end{document}